\documentclass[letterpaper, 10 pt, conference]{ieeeconf}

\IEEEoverridecommandlockouts 

\usepackage{stackrel}
\usepackage{graphicx}
\usepackage{siunitx}
\graphicspath{{./Images/}}
\usepackage{amssymb,amsmath}
\usepackage{multicol}
\usepackage{xcolor}
\usepackage{bm} 
\newtheorem{remark}{Remark}
\newtheorem{definition}{Definition}
\newtheorem{lemma}{Lemma}
\newtheorem{theorem}{Theorem}

\newtheorem{proofoflemma}{Proof of Lemma}
\newtheorem{proofoftheorem}{Proof of Theorem}
\newcommand{\rd}[1]{{\color{black}{#1}}}

\title{\LARGE \bf A Dynamic Observation Strategy for \\ Multi-agent Multi-armed Bandit Problem}

\author{Udari Madhushani
and Naomi Ehrich Leonard
	\thanks{This research has been supported in part by ARO grant W911NF-18-1-0325 and ONR grant N00014-19-1-2556. Department of Mechanical and Aerospace Engineering, Princeton University, NJ 08544, USA.
		{\tt\small \{udarim,naomi\}@princeton.edu}}%
}

\begin{document}

\maketitle
\thispagestyle{empty}
\pagestyle{empty}

\begin{abstract}
We define and analyze a multi-agent multi-armed bandit problem in which decision-making agents can observe the choices and rewards of their neighbors under a linear observation cost.  Neighbors are defined by a network graph that encodes the inherent observation constraints of the system. We define a cost associated with observations such that at every instance an agent makes an observation it receives a constant observation regret. We design a sampling algorithm and an observation protocol for each agent to maximize its own expected cumulative reward through minimizing expected cumulative sampling regret and expected cumulative observation regret. For our proposed protocol, we prove that total cumulative regret is logarithmically bounded. We verify the accuracy of analytical bounds using numerical simulations. 
\end{abstract}

\allowdisplaybreaks

\section{Introduction} \label{sect:Introduction}
The effect of communication structure in cooperative and competitive multi-agent systems has been extensively studied in decision theory. Performance of a group of social learners can be improved by the shared information among individuals. In most real-world decision-making processes, however, information sharing between agents can be costly. As a result, directed communication, where each agent only needs to observe its neighbors, has advantages over undirected communication, where each agent sends and receives information. Even when observation costs are high, agents can keep costs to a minimum by choosing when and whom to observe as a function of their own performance. Further, in this setting costs associated with cooperation can be avoided.
    
Consider the problem of a group of fishermen foraging in an uncertain environment that consists of a distribution of spatial resource (fish). Because of the natural dynamics of fish, environmental conditions, and other external factors, the resource will be distributed stochastically. As a result, a fisherman will receive different reward values (number of fish harvested) at different times, even when sampling from the same patch. 
Thus, in order to maximize cumulative reward fishermen need to be able to exploit, i.e., forage in well sampled patches known to provide better harvest,
and to explore, i.e, forage in poorly sampled patches, which is riskier 
but may provide even better harvest than well sampled patches. Benefiting from exploitation requires sufficient exploration and identification of the patches that yield highest rewards. 
More generally, optimal foraging performance comes from balancing the trade-off between exploring and exploiting. This is known as the explore-exploit dilemma. 

Multi-armed bandit (MAB) problems are a set of mathematical models that have been proposed to capture the salient features of explore-exploit trade-offs \cite{sutton1998introduction,robbins1952some}. For the standard MAB problem the reward distributions associated with options are static. An agent estimates the expected reward of each option using the rewards it receives through sampling. The agent chooses among options by considering a trade-off between estimated expected reward (exploiting) and the uncertainty associated with the estimate (exploring). Therefore, in the frequentist setting, the natural way of estimating the expectation of the reward is to consider the sample average \cite{lai1985asymptotically,agrawal1995sample,auer2002finite}. The papers \cite{kaufmann2012bayesian,reverdy2014modeling} present how to incorporate prior knowledge about reward expectation in the estimation step by leveraging the theory of conditional expectation in the Bayesian setting. 

Multi-agent multi-armed bandit (MAMAB) problems consider a group of individuals facing the same MAB problem simultaneously. For an individual to maximize its own reward, it will naturally seek to observe its neighbors and use those observations to improve its performance. Individual and group performance of agents will vary according to the observation structure, i.e., who is observing whom, and the type of information they observe. For example, if the agents are cooperative and can broadcast signals, they could share their estimates of rewards.  
When there are constraints, such as communication costs and privacy concerns,
they might instead share only their instantaneous rewards and choices.  Even without the ability to broadcast, agents may still be able to use sensors to observe the instantaneous rewards and choices of neighbors. 
A centralized multi-agent setting is considered in  \cite{anantharam1987asymptotically} and a decentralized setting is considered in \cite{kalathil2014decentralized}. The papers \cite{landgren2016distributed,landgren2016distributedCDC} use a running consensus algorithm in which agents observe the reward estimates of their neighbors. In \cite{kolla2018collaborative,landgren2018social} an MAMAB problem is studied in which agents observe instantaneous rewards and choices in a leader-follower setting.

In all of these previous works, communication between agents is assumed to be cost free. However, in real world settings observing neighbors or exchanging information with neighbors is costly. In the present paper, we propose a setting in which agents can decide when and whom to observe in order to receive maximum benefits from observations that incur a cost. An underlying undirected network graph defines neighbors and models the inherent observation constraints present in the network.  Agents receive a fixed observation cost at every instance they observe a neighbor. 

To account for the observation cost, we define cumulative regret to be the total cumulative regret agents receive from sampling suboptimal options (sampling regret) and from observing neighbors (observation regret).  Deterministic \cite{landgren2018social} and probabilistic \cite{madhushani2019heterogeneous} communication strategies proposed in the MAB literature lead to a linear cumulative observation regret.  Our main contribution is the design of a new strategy for which we prove a logarithmic total cumulative regret, i.e., order-optimal performance.  Our design leverages the intuition that it is most useful to observe neighbors when uncertainty associated with estimations of rewards is high. 

In Section \ref{Secn:MAMAB} we introduce the MAMAB problem 
and we propose an efficient sampling rule and a communication protocol for an agent to maximize its own total expected cumulative reward. We analyze the performance of the proposed sampling rule in Section \ref{Secn:Regret}. In Section \ref{SubSecn:Regret} we analytically upper bound the expected cumulative regret and in Section \ref{SubSecn:CObsRegret} we analytically upper bound the expected observation regret. We present the upper bound for the total expected cumulative regret in section \ref{SubSec:TotalRegret}. In Section \ref{Secn:Simu} we provide numerical simulation results and computationally validate the analytical results. 
We conclude in Section \ref{Secn:Concl} and provide additional mathematical details in the Appendix.


\section{Multi-agent Multi-armed Bandit Problem}\label{Secn:MAMAB}

In this section we present the mathematical formulation of the MAMAB problem studied here. Let $N$ be the number of options (arms) and $K$ the number of agents. Define $X_i$ as the random variable that denotes reward associated with option $i\in\mathcal{I}=\{1,2,\ldots,N\}$. In this paper we assume that all the reward distributions are sub-Gaussian. Let $\sigma_i$ be the variance proxy of $X_i,$ 
and $\mu_i$ the expected reward of option $i.$ Let $i^*$ be the optimal option with highest expected reward $\mu_{i^*}=\max\{\mu_1,\mu_2,\ldots,\mu_N\}.$ Each agent $k \in \{1,\ldots, K\}$ chooses one option at each time step $t \in \{1,2,\ldots,T\}$ with the goal of minimizing its cumulative regret. In MAB problems, cumulative regret is typically defined as cumulative sampling regret, which is equivalent to  expected number of times suboptimal options are selected.  We let cumulative regret  be the sum of cumulative sampling regret and a cumulative observation regret that accumulates a fixed cost for every observation of a neighbor.

We assume that the expected reward values $\mu_{i}$  are unknown and the variance proxy values $\sigma_{i}$ are known to the agents. To improve its own performance, each agent observes its neighbors according to an observation protocol that we define. We use a network graph to encode hard observation constraints and this defines neighbors of agents. Let $\mathcal{G}(\mathcal{V},\mathcal{E})$ be an undirected graph. $\mathcal{V}$ is a set of $K$ nodes, such that node $k$ in $\mathcal{V}$ corresponds to agent $k$ for $k \in \{1, \ldots, K\}$.  $\mathcal{E}$ is a set of edges between nodes in $\mathcal{V}$. 
If there is an edge $e(k,j) \in\mathcal{E}$ between node $k$ and node $j$, then we say that agent $k$  and agent $j$ are neighbors. Since the graph is undirected, $e(k,j)\in\mathcal{E}\iff e(j,k)\in\mathcal{E}.$
Let $d_k$ be the number of neighbors of agent $k.$


Let $\varphi_{k}^{t}\in \mathcal{I}$ and $X_{k}^{t}$ be random variables that denote the option chosen by   agent $k$ and the reward received by  agent $k$ at time $t$, respectively.
 Let $\mathbb{I}_{\{\varphi_{k}^{t} =i\}}$ be a random variable that takes value 1 
 if  option $i$ is chosen by agent $k$ at time $t$ and is 0  otherwise. Let $\mathbb{I}^{t}_{\{k,j\}}$ be a
random variable that takes value 1
 if agent $k$ can observe agent $j$ at time $t$ and is 0 otherwise.

 In order to maximize the cumulative reward in the long run, agents need to both identify the best options through exploring and sample the best options through exploiting. Observing neighbors allows an agent to receive more information about options and hence obtain better estimates about expected reward values of options. This leads to less exploring and more exploiting, which reduces the regret an agent receives due to sampling suboptimal options. However, since taking observations is costly, an agent is required to find a trade-off between the information gain and the cost associated with observations. Let $c_{k,j}$ be the cost incurred by agent $k$ when it observes the instantaneous reward and choice of agent $j$ at time step $t$. In this paper we consider the case in which $c_{k,j}=c, \forall j,k.$

Let the number of times that agent $k$ samples option $i$ until time $t$ be given by the
random variable  $n_{i}^{k}(t) = \sum_{\tau=1}^t\mathbb{I}_{\{\varphi_{k}^{\tau} =i\}}$. And let the total number of times that  agent $k$ observes rewards from option $i$ until time $t$ be given by the
random variable $N_{i}^{k}(t)$, where
\begin{align*}
N_{i}^{k}(t)=\sum_{\tau=1}^{t}\sum_{j=1}^{K}\mathbb{I}_{\{\varphi_{j}^{\tau} =i\}}\mathbb{I}^{\tau}_{\{k,j\}}.
\end{align*}
 
We define a sampling rule based on the well known UCB (Upper Confidence Bound) rule for a single agent \cite{auer2002finite}.  The UCB rule chooses the option at time $t$ that maximizes an objective function that is the sum of an exploit term, equal to the estimate of the reward mean at time $t$, and an explore term, equal to a measure of uncertainty in that estimate at time $t$. Our sampling rule for agent $k$ in the MAMAB problem accounts for the observations of neighbors by using them to improve its estimate and reduce its uncertainty. Let the
estimate by agent $k$ of the expected reward from option $i$  at time $t$ be given by the random variable $\widehat{\mu}_{i}^{k}(t)$, where
\begin{align*}
\widehat{\mu}_{i}^{k}(t)=\frac{S_{i}^{k}(t)}{N_{i}^{k}(t)},
\end{align*}
and $S_{i}^{k}(t)=\sum_{\tau=1}^{t}\sum_{j=1}^{K}X_{i}\mathbb{I}_{\{\varphi_{j}^{\tau} =i\}}\mathbb{I}^{\tau}_{\{k,j\}}$ is the total reward observed by agent $k$ from option $i$ until time $t$.

\begin{definition}\label{def:UCB based}
The {\em sampling rule} $\{\varphi_{k}^{t}\}_1^{T}$ {\em for agent $k$ at time} $t \in \{1, \ldots, T\}$ is defined as
	\begin{align}
	\mathbb{I}_{\{\varphi_{k}^{t+1}=i\}}=\left\{
	\begin{array}{cl} 1 &, \:\:\:Q^{k}_i(t)=\max\{Q^{k}_1(t),\cdots,Q^{k}_N(t)\}\label{eq:UCBallocation}\\
  	0 &, \:\:\: {\mathrm{o.w.}}\end{array}\right.
	\end{align}
	with 	
	\begin{align}
	Q^{k}_i(t)&= \widehat{\mu}_{i}^{k}(t)+C_i^k(t)\label{eq:UCBQ}\\
	C_i^k(t)&=\sigma_{i}\sqrt{2(\xi+1)\frac{\log t}{N_i^k(t)}}, \label{eq:Uncertainity}
	\end{align}
where $\xi>1$ is a tuning parameter that captures the trade-off between exploring and exploiting.
\end{definition}

To find a balance between information gain and observation cost we define an observation rule for agents so that they  choose to incur the cost of making observations of neighbors only when observations are most needed, i.e., when their own uncertainty is high.
In the following observation rule, an agent observes the instantaneous rewards and choices of all of its neighbors only when it is exploring, since it explores when uncertainty is high.  If  agent $k$ chooses the option at time $t$ that corresponds to the maximum of its estimates of reward means, $\widehat{\mu}_1^k(t), \ldots, \widehat{\mu}_N^k(t)$, then it is exploiting and it does not observe its neighbors.
\begin{definition}\label{def:observation rule}
The {\em observing rule} $\mathbb{I}_{\{k,j\}}^t$ {\em for agent $k$ at time} $t \in \{1, \ldots, T\}$ and $\forall j$ is defined as
	\begin{align}
	\mathbb{I}_{\{k,j\}}^{t+1}=\left\{
	\begin{array}{cl} 0\!\!\!\!  &, \varphi_k^t=i,\: s.t.\: \widehat{\mu}_i^k(t)\!=\!\max\{\widehat{\mu}_{1}^{k}(t),\cdots,\widehat{\mu}_{N}^{k}(t)\}\label{eq:UCBobservation}\\
  	1\!\!\!\!  &, {\mathrm{o.w.}}\end{array}\right.
	\end{align}
\end{definition}

\section{Performance Analysis}\label{Secn:Regret}
In this section we analyze the cumulative regret of agent $k$ due to sampling suboptimal options and observing neighbors when employing the sampling rule of Definition~\ref{def:UCB based} and observation rule of Definition~\ref{def:observation rule}.

\subsection{Sampling Regret Analysis}\label{SubSecn:Regret}
Let $i$ be a suboptimal option. The total number of times agent $k$ samples from  option $i$ can be upper bounded as 
\[
n_{i}^{k}(T)=\sum_{t=1}^T\mathbb{I}_{\{\varphi_{k}^{t} =i\}}
\leq \sum_{t=1}^{T}\mathbb{I}_{\{Q_{i}^{k}(t)\geq Q_{i^{*}}^{k}(t)\}}.
\]
Here $\mathbb{I}_{\{Q_{i}^{k}(t)>Q_{i^{*}}^{k}(t)\}}$ is an indicator function such that 
\begin{align*}
\mathbb{I}_{\{Q_{i}^{k}(t)>Q_{i^{*}}^{k}(t)\}}=\left\{
	\begin{array}{cl} 1 &, \:\:\:Q^{k}_i(t)\geq Q^{k}_{i^{*}}(t)\\
  	0 &, \:\:\: {\mathrm{o.w.}}\end{array}\right.
\end{align*}
Thus we have 
\begin{align*}
\mathbb{E}\left(n_{i}^{k}(T)\right)\leq \sum_{t=1}^{T}\mathbb{P}\left(Q_{i}^{k}(t)\geq Q_{i^{*}}^{k}(t)\right). 
\end{align*}

Let $R_{s}^{k}(T)$ be the cumulative sampling regret of agent $k$ from option $i$ until time $T$. Recall that the cumulative regret is defined as the loss incurred by sampling suboptimal options. Define $\Delta_{i}= \mu_{i^{*}}-\mu_{i}.$ Then we have, from \cite{lai1985asymptotically},
\begin{align}
\mathbb{E}\left(R_{s}^{k}(T)\right)=\sum_{i=1}^N\Delta_{i}\mathbb{E}\left(n_{i}^{k}(T)\right).\label{eq:regretValue}
\end{align}
To analyze the expected number of samples from suboptimal options until time $T$, we first note that $\forall i,k,t$ we have 
\begin{align*}
& \left\{\mathbb{I}_{\{\varphi_k^{t+1}=i\}}\right\}\subseteq \left\{Q_{i}^{k}(t)\geq Q_{i^{*}}^{k}(t)\right\}\subseteq \left\{\mu_{i^{*}}<\mu_{i}+2C_{i}^{k}(t)\right\}\nonumber\\
&\:\:\:\: \:\:\:\:\:\:\:\:\cup\left\{\widehat{\mu}_{i^{*}}^{k}(t)\leq \mu_{i^{*}}-C_{i^{*}}^{k}(t)\right\}\cup \left\{\widehat{\mu}_{i}^{k}(t)\geq \mu_{i}+C_{i}^{k}(t)\right\}
\end{align*}
and so
\begin{align}
&\mathbb{E}\left(n_{i}^{k}(T)\right) \leq \sum_{t=1}^{T}\mathbb{P}\left(\mu_{i^{*}}<\mu_{i}+2C_{i}^{k}(t)\right) +\nonumber\\
&\sum_{t=1}^{T}\mathbb{P}\!\left(\widehat{\mu}_{i^{*}}^{k}(t)\leq \mu_{i^{*}}-C_{i^{*}}^{k}(t)\right) \! + \!\! \sum_{t=1}^{T}\mathbb{P}\!\left(\widehat{\mu}_{i}^{k}(t)\geq \mu_{i}+C_{i}^{k}(t)\right)
.\label{eq:Prob}
\end{align}

Next we analyze concentration probability bounds on the estimates of options.

\begin{theorem}\label{thm:TailProb}
For any $\zeta>1$ and for $\sigma_i>0$ there exists a $\vartheta>0$ such that
\begin{align*}
    \mathbb{P}\left(\widehat{\mu}_i^k(T)-{\mu}_i>\sqrt{\frac{\vartheta}{N_{i}^{k}(T)}}\right)\leq \frac{\nu\log (d_k+1)T}{\exp(2\kappa \vartheta)}
\end{align*}
where 
\begin{align*}
\nu=\frac{1}{\log \zeta},\:\:\:\: \kappa=\frac{1}{\sigma_i^2\left(\zeta^{\frac{1}{4}}+\zeta^{-\frac{1}{4}}\right)^2}.
\end{align*}
\end{theorem}

The proof of Theorem~\ref{thm:TailProb} can be found in the paper \cite{madhushani2019heterogeneous}. Using symmetry we conclude that
\begin{align*}
    \mathbb{P}\left(\Big |\widehat{\mu}_i^k(T)-{\mu}_i\Big |>\sqrt{\frac{\vartheta}{N_{i}^{k}(T)}}\right)\leq \frac{\nu\log (d_k+1)T}{\exp(2\kappa \vartheta)}.
\end{align*}

\begin{lemma}\label{lem:TailProb}
For $\vartheta=2\sigma_{i}^{2}(\xi+1)\log T$ and $\xi>1$ there exists a $\zeta>1$ such that
\begin{align*}
    \mathbb{P}\left(\Big|\widehat{\mu}_i^k(T)-{\mu}_i\Big |>\sigma_{i}\sqrt{\frac{2(\xi+1)\log T}{N_{i}^{k}(T)}}\right) \! \leq \! \frac{\nu\log (d_k+1)T}{T^{\xi+1}}.
\end{align*}
\end{lemma}
The proof of Lemma~\ref{lem:TailProb} can be found in the paper \cite{madhushani2019heterogeneous}.

We proceed to upper bound the summation of the probabilities of the events $\left\{\mu_{i^{*}}<\mu_{i}+2C_{i}^{k}(t)\right\}$ for $t\in\{1,2,\ldots,T\}$ as follows. Using equation (\ref{eq:Uncertainity}) we have that the inequality
$
\mu_{i^{*}}<\mu_{i}+2C_{i}^{k}(t)
$
implies
\begin{align*}
\frac{\Delta_i^2}{4\sigma_{i}^2}\left(N_i^k(t)\right)^2-2(\xi+1)\log t \left(N_i^k(t)\right)&<0.
\end{align*}
This inequality does not hold for $N_{i}^{k}(t)> \eta_i(t)$, where
\begin{align*}
\eta_i(t) =\frac{8\sigma_{i}^{2}(\xi+1)}{\Delta_{i}^{2}}\log t
.
\end{align*}
Thus we have
\begin{align}
\sum_{t=1}^{T}\mathbb{P}\left(Q_{i}^{k}(t)\geq Q_{i^{*}}^{k}(t),N_{i}^{k}(t)>\eta_i(t)\right)\leq \eta_i(T). \label{eq:eta}
\end{align}

From the probability bounds given in Lemma~\ref{lem:TailProb} and  (\ref{eq:eta}), the total expected number of times agent $k$ samples suboptimal option $i$ until  time $T$ is upper bounded as
\begin{align}
\mathbb{E}\left(n_{i}^{k}(T)\right)\leq&  \frac{1}{\log \zeta}(1+\log (d_k+1))+\frac{8\sigma_{i}^{2}(\xi+1)}{\Delta_{i}^{2}} \log T\nonumber\\
&+\frac{1}{2^{\xi}\log \zeta}\left(\frac{\log (d_k+1)}{\xi}+\frac{2}{\xi-1}\right)\nonumber\\
&+\frac{1}{T^{\xi-1}\log \zeta}\left(\frac{\log (d_k+1)}{T\xi}+\frac{1}{\xi-1}\right)\label{eq:suboptions}
\end{align}
where $\zeta,\xi>1.$

From equation (\ref{eq:regretValue}) the expected cumulative sampling regret of agent $k$ until time $T$ is upper bounded as
\begin{align}
\mathbb{E}\left(R_{s}^k(T)\right)\leq&  \sum_{i=1}^N\frac{\Delta_i}{\log \zeta}(1+\log (d_k+1))\nonumber\\
&+\frac{8\sigma_{i}^{2}(\xi+1)}{\Delta_{i}} \log T\nonumber\\
&+\sum_{i=1}^N\frac{\Delta_i}{2^{\xi}\log \zeta}\left(\frac{\log (d_k+1)}{\xi}+\frac{2}{\xi-1}\right)\nonumber\\
&+\sum_{i=1}^N\frac{\Delta_i}{T^{\xi-1}\log \zeta}\left(\frac{\log (d_k+1)}{T\xi}+\frac{1}{\xi-1}\right).
\end{align}

\subsection{Observation Regret Analysis}\label{SubSecn:CObsRegret}

Recall that $c$ is the constant unit cost associated with observations. Let $R_{o}^{k}(T)$ be the cumulative observation regret of agent $k$ at time step $T.$ Then we have
\begin{align*}
R_{o}^k(T)=c\sum_{t=1}^T\sum_{j=1}^K\mathbb{I}^t_{\{k,j\}}.
\end{align*}
This is equivalent to the number of observations taken by agent $k$ until time $T.$  Expected cumulative observation regret can be expressed as
\begin{align}
\mathbb{E}\left(R_{o}^k(T)\right)=c\sum_{t=1}^T\sum_{j=1}^K\mathbb{E}\left(\mathbb{I}^t_{\{k,j\}}\right).\label{eq:ExpectCumObsRegret}
\end{align}
So expected cumulative observation regret can be upper bounded by upper bounding the expected number of observations until time $T$:
\begin{align}
\sum_{t=1}^T&\sum_{j=1}^K\mathbb{E}\left(\mathbb{I}_{\{k,j\}}\right)\nonumber \\
&=d_k\sum_{t=1}^T\mathbb{P}\left(\varphi_k^t=i,\widehat{\mu}_i^k(t)\neq \max\{\widehat{\mu}_{1}^{k}(t),\cdots,\widehat{\mu}_{N}^{k}(t)\}\right).\label{eq:ExpectNumObs}
\end{align}

To analyze the expected number of observation, we use
\begin{align*}
    \mathbb{P}&\left(\varphi_k^t=i,\widehat{\mu}_i^{k}(t)\neq \max\{\widehat{\mu}_{1}^{k}(t),\cdots,\widehat{\mu}_{N}^{k}(t)\}\right)=\\
    & \mathbb{P}\left(\varphi_k^t=i^*,\widehat{\mu}_{i^*}^{k}(t)\neq \max\{\widehat{\mu}_{1}^{k}(t),\cdots,\widehat{\mu}_{N}^{k}(t)\}\right)\\
    &+ \mathbb{P}\left(\varphi_k^t= i,\widehat{\mu}_i^{k}(t)\neq \max\{\widehat{\mu}_{1}^{k}(t),\cdots,\widehat{\mu}_{N}^{k}(t)\},i\neq i^*\right).
\end{align*}
We first upper bound the expected number of times agent $k$ observes its neighbors until time $T$ when it decides to explore after sampling a suboptimal option.
\begin{lemma}\label{lem:ExploreSuboptimal}
For all suboptimal $i\neq i^*$ we have
\begin{align*}
&\sum_{t=1}^T\mathbb{P}\left(\varphi_k^t= i,\widehat{\mu}_i^{k}(t)\neq \max\{\widehat{\mu}_{1}^{k}(t),\cdots,\widehat{\mu}_{N}^{k}(t)\},i\neq i^*\right)\\
&\leq \frac{N-1}{\log \zeta}(1+\log (d_k+1))+\sum_{\stackrel{i=1}{i\neq i^*}}^N\frac{8\sigma_{i}^{2}(\xi+1)}{\Delta_{i}^{2}} \log T\\
&+\frac{N-1}{2^{\xi}\log \zeta}\left(\frac{\log (d_k+1)}{\xi}+\frac{2}{\xi-1}\right)\\
&+\frac{N-1}{T^{\xi-1}\log \zeta}\left(\frac{\log (d_k+1)}{T\xi}+\frac{1}{\xi-1}\right).
\end{align*}
\end{lemma}

The proof of Lemma \ref{lem:ExploreSuboptimal} is given in the Appendix.

Next we analyze the expected number of times agent $k$ observes its neighbors until time $T$ when it decides to explore after sampling the optimal option.

Note that $\forall i,k,t$ we have
\begin{align*}
   & \{\varphi_k^t=i^*,\widehat{\mu}_{i^*}\neq \max\{\widehat{\mu}_{i}^{k}(t),\cdots,\widehat{\mu}_{N}^{k}(t)\}\}\subseteq \\
    &\{\widehat{\mu}_{i^*}^{k}(t)\leq \mu_{i^*}-C_{i^*}^{k}(t)\}\\
    &\cup\{\widehat{\mu}_{i^*}^{k}(t)\geq \mu_{i^*}-C_{i^*}^{k}(t),\exists i, s.t. (\widehat{\mu}_{i}^k(t)\geq \mu_{i^*}-C_{i^*}^{k}(t)\}.
\end{align*}
Thus we have
\begin{align*}
    &\sum_{i=1}^T\mathbb{P}\left(\varphi_k^t=i^*,\widehat{\mu}_{i^*}\neq \max\{\widehat{\mu}_{i}^{k}(t),\cdots,\widehat{\mu}_{N}^{k}(t)\}\right)\\
    &\leq \sum_{i=1}^T\mathbb{P}\left(\widehat{\mu}_{i^*}^{k}(t)\leq \mu_{i^*}-C_{i^*}^{k}(t)\right) +\\
    &\sum_{i=1}^T\! \mathbb{P}\!\left(\widehat{\mu}_{i^*}^{k}(t)\geq \mu_{i^*}-C_{i^*}^{k}(t),\exists i, s.t. (\widehat{\mu}_{i}^k(t)\geq \widehat{\mu}_{i^*}^{k}(t)\right).
\end{align*}
From Lemma \ref{lem:TailProb} we have
\begin{align}
\sum_{i=1}^T&\mathbb{P}\left(\widehat{\mu}_{i^*}^{k}(t)\leq \mu_{i^*}-C_{i^*}^{k}(t)\right)\leq 
\frac{1}{\log \zeta}(1+\log (d_k+1))\nonumber\\
&+\frac{1}{2^{\xi}\log \zeta}\left(\frac{\log (d_k+1)}{\xi}+\frac{2}{\xi-1}\right)\nonumber\\
&+\frac{1}{T^{\xi-1}\log \zeta}\left(\frac{\log (d_k+1)}{T\xi}+\frac{1}{\xi-1}\right).\label{eq:TailBound}
\end{align}

\begin{theorem}\label{thm:explSubO} For all suboptimal options $i\neq i^*$ we have
\begin{align*}
    \sum_{i=1}^T&\mathbb{P}\left(\widehat{\mu}_{i^*}^{k}(t)\geq \mu_{i^*}-C_{i^*}^{k}(t),\exists i, s.t. (\widehat{\mu}_{i}^k(t)\geq \widehat{\mu}_{i^*}^{k}(t)\right)\leq\\
    &\sum_{i=1}^N\frac{8\sigma_{i}(\xi+1)}{\Delta^2_{i}}\log T +\frac{N-1}{2^{\xi}\log \zeta}\left(\frac{\log (d_k+1)}{\xi}+\frac{2}{\xi-1}\right)\\
    &+\frac{N-1}{\log \zeta}(1+\log (d_k+1))\\
&+\frac{N-1}{T^{\xi-1}\log \zeta}\left(\frac{\log (d_k+1)}{T\xi}+\frac{1}{\xi-1}\right).
\end{align*}
\end{theorem}
The proof of Theorem~\ref{thm:explSubO} is given in the Appendix.

Now we proceed to state the main result of this paper, which is that the total expected cumulative observation regret until time $T$ for agent $k$ employing the sampling rule given by Definition~\ref{def:UCB based} and the observation rule given by Definition~\ref{def:observation rule} is upper bounded logarithmically in $T$. 
\begin{theorem}\label{thm:ReglSubO} Expected cumulative observation regret until time $T$ for agent $k$ can be upper bounded as
\begin{align*}
\mathbb{E}&\left(R_o^k(T)\right)\leq \sum_{i=1}^N\frac{8\sigma_{i}(\xi+1)}{\Delta^2_{i}}\log T\\
&+\frac{cd_k(2N-1)}{\log \zeta}(1+\log (d_k+1))\\
&+\frac{cd_k(2N-1)}{2^{\xi}\log \zeta}\left(\frac{\log (d_k+1)}{\xi}+\frac{2}{\xi-1}\right)\\
&+\frac{cd_k(2N-1)}{T^{\xi-1}\log \zeta}\left(\frac{\log (d_k+1)}{T\xi}+\frac{1}{\xi-1}\right).
\end{align*}
\end{theorem}

Theorem \ref{thm:ReglSubO} follows from equations (\ref{eq:ExpectCumObsRegret})-(\ref{eq:TailBound}), Lemma \ref{lem:ExploreSuboptimal} and Theorem \ref{thm:explSubO}.

\begin{remark}\label{Rem:LinearRegret}
Note that for deterministic communication strategies \cite{landgren2016distributed,landgren2018social} the expected cumulative observation regret until time $T$ for agent $k$ is linear in $T$:
\begin{align*}
\mathbb{E}\left(R_{o}^k(T)\right)=c\sum_{t=1}^T\sum_{j=1}^K\mathbb{E}\left(\mathbb{I}^t_{\{k,j\}}\right)=cd_kT.
\end{align*}
For the probabilistic observation strategy of \cite{madhushani2019heterogeneous} the expected cumulative observation regret until time $T$ for agent $k$ is linear in $T$:
\begin{align*}
\mathbb{E}\left(R_{o}^k(T)\right)=c\sum_{t=1}^T\sum_{j=1}^K\mathbb{E}\left(\mathbb{I}^t_{\{k,j\}}\right)=cd_kp_kT,
\end{align*}
where $p_k$ is the observation probability of agent $k.$ Thus, our proposed sampling rule and observation rule outperform these strategies when there are cumulative observation costs.
\end{remark}

\subsection{Total expected cumulative regret}\label{SubSec:TotalRegret}
Total expected cumulative regret $\mathbb{E}\left(R^k(T)\right)$  is defined as the summation of expected cumulative sampling regret and expected cumulative observation regret until time $T$:
\begin{align*}
 \mathbb{E}\left(R^k(T)\right)  = \sum_{i=1}^N\mathbb{E}\left(R_{i}^{k}(T)\right)+\mathbb{E}\left(R_o^k(T)\right).
\end{align*}
Let $\sum_{i=1}^N\Delta_i=\widetilde{\Delta}.$ Total expected cumulative regret until time $T$ of agent $k$ is upper bounded as
\begin{align}
& \mathbb{E}\left(R_{s}^k(T)\right)\leq \sum_{\stackrel{i=1}{i\neq i^*}}^N\frac{8\sigma_{i}^{2}(\xi+1)}{\Delta^2_{i}} \log T \nonumber \\&  \frac{\widetilde{\Delta}+cd_k(2N-1)}{\log \zeta}(1+\log (d_k+1))\nonumber\\
&+\frac{\widetilde{\Delta}+cd_k(2N-1)}{2^{\xi}\log \zeta}\left(\frac{\log (d_k+1)}{\xi}+\frac{2}{\xi-1}\right)\nonumber\\
&+\frac{\widetilde{\Delta}+cd_k(2N-1)}{T^{\xi-1}\log \zeta}\left(\frac{\log (d_k+1)}{T\xi}+\frac{1}{\xi-1}\right).
\end{align}

\section{Simulation Results}\label{Secn:Simu}
In this section we present numerical simulation results for a network of 6 agents with underlying observation structure defined by the star graph: {\rd{the center agent observes all other agents and all other agents only observe the center agent. Agents other than the center agent are interchangeable and their average regret and individual regret are the same. }} We present numerical simulations \rd{to evaluate} the performance of the sampling rule and observation rule given by Definitions~\ref{def:UCB based} and \ref{def:observation rule}. 

The 6 agents play the same MAB problem with 10 options. In all simulations  the reward distributions are Gaussian with variance $\sigma_i = 5$, $i=1, \ldots, 10$, and mean values:
\begin{center}
\begin{tabular}{*{11}{|c}|}
\hline
i & 1 & 2 & 3 & 4 & 5 & 6 & 7 & 8 & 9 & 10\\
\hline
$\mu_i$ & 40 & 50 & 50 & 60 & 70 & 70 & 80 & 90 & 92 & 95 \\
\hline
\end{tabular}.
\end{center} 
The communication cost $c = 1$. We set the sampling rule parameter $\xi= 1.01$. We provide results for 1000 time steps with 1000 Monte Carlo simulations. 

{\rd{
Figure \ref{Fig:ObsRule} shows simulation results for the expected cumulative sampling regret of a group of 6 agents using the proposed sampling and observation rules. The blue dashed line shows regret of the center agent. The green dash-dot line shows the average regret of the agents not in the center. The red dotted line shows the average expected cumulative sampling regret over all agents. It can be observed that the expected cumulative sampling regret is logarithmic in time.}} For comparison, we plot the average expected cumulative regret of the agents when they make no observations of neighbors (solid gold line).  
When agents are not making observations they are interchangeable, and so 
the average performance and the individual performance are the same. {\rd{The simulation results illustrate that the performance of every agent improves significantly when it observes neighbors according to the proposed protocol.}} The simulation results further show that the center agent outperforms the other agents. This is to be expected since the center agent has more neighbors than the other agents.

{\rd{Figure \ref{Fig:NoObs} shows simulation results for expected observation regret. 
It can be seen that the expected observation regret is logarithmic in time, as proved in Theorem \ref{thm:ReglSubO}. Since the center agent has more neighbors than the others agents, its observation regret is the highest. However, the results illustrate that when observation cost is small, a significant performance improvement can be obtained for a small observation regret.
\begin{figure}[t]
    \centering
    \includegraphics[width=0.4\textwidth]{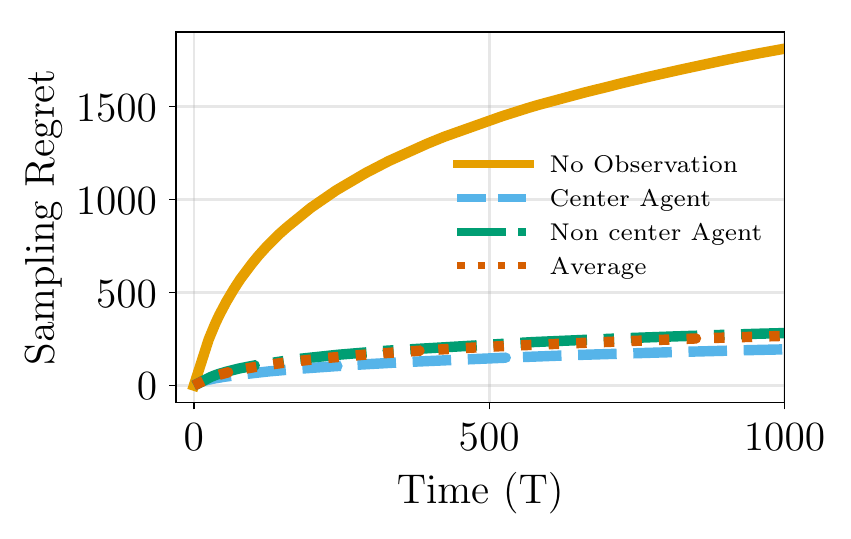}
    \caption{Dashed and dotted lines show expected  cumulative sampling regret of the agents using the sampling rule and observation rule of Definitions~\ref{def:UCB based} and \ref{def:observation rule}
    with underlying star observation structure.  The solid line shows the average performance of agents when they are not observing their neighbors.}
\label{Fig:ObsRule}
\end{figure}

\begin{figure}[t]
    \centering
    \includegraphics[width=0.4\textwidth]{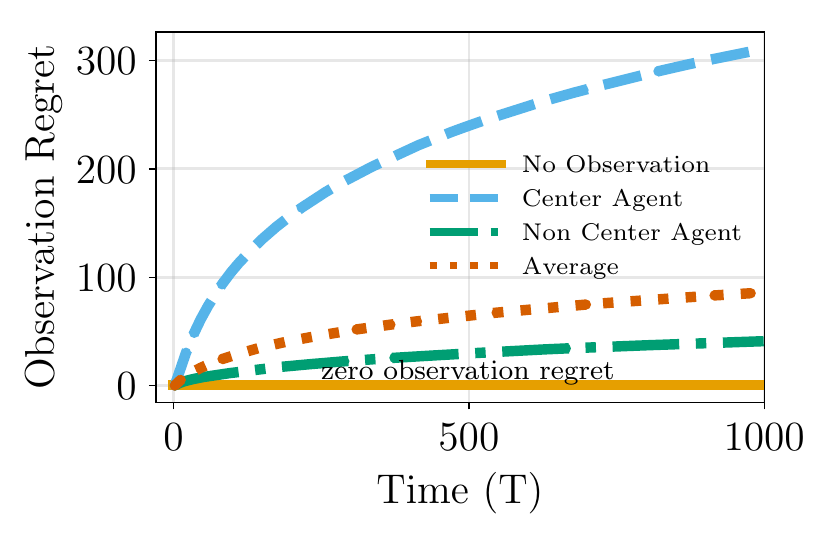}
    \caption{Dashed and dotted lines show expected  cumulative observation regret of the agents using the sampling rule and observation rule of Definitions~\ref{def:UCB based} and \ref{def:observation rule}
    with underlying star observation structure. The solid line shows that agents do not suffer from any observation regret when they do not observe their neighbors.} 
\label{Fig:NoObs}
\end{figure}
}}


\section{Conclusions}\label{Secn:Concl}
We studied an MAMAB problem where agents can observe the instantaneous choices and rewards of their neighbors but incur a cumulative cost each time they make an observation of a neighbor.  We proposed a sampling rule and an observation rule in which an agent observes its neighbors only when it has decided to explore.  We defined total expected cumulative regret to be the regret agents receive due to sampling suboptimal options and to observing neighbors. Deterministic and stochastic observation strategies for MAB protocols in the literature yield an expected cumulative observation regret that is linear in time $T$. We analytically proved that under the proposed sampling and observation rules, expected cumulative regret of each agent is bounded logarithmically in $T$. Accuracy of the upper bound has been verified computationally through numerical simulations.


\begin{appendix}

\setcounter{proofoflemma}{1}
\begin{proofoflemma}
Note that $\forall i,k,t$ we have
\begin{align*}
 \mathbb{P}\left(\varphi_k^t= i,\widehat{\mu}_i^{k}(t)\neq \max\{\widehat{\mu}_{1}^{k}(t),\cdots,\widehat{\mu}_{N}^{k}(t)\},i\neq i^*\right)\leq \\
  \mathbb{E}\left(\mathbb{I}_{\{\varphi_k^t= i\}}\right).
\end{align*}
Then we have
\begin{align*}
    \sum_{t=1}^T\mathbb{P}\left(\varphi_k^t= i,\widehat{\mu}_i^{k}(t)\neq \max\{\widehat{\mu}_{1}^{k}(t),\cdots,\widehat{\mu}_{N}^{k}(t)\},i\neq i^*\right)\\
    \leq \sum_{t=1}^T\sum_{i=1}^N\mathbb{E}\left(\mathbb{I}_{\{\varphi_k^t= i\}}\right).
\end{align*}
Lemma \ref{lem:ExploreSuboptimal} follows from equation (\ref{eq:suboptions}).
\end{proofoflemma}

\setcounter{proofoftheorem}{1}
\begin{proofoftheorem}
Let $i$ be a suboptimal option with highest estimated expected reward for agents $k$ at time $t.$ Then we have $i=\arg\max \{\widehat{\mu}_{1}^{k}(t),\cdots,\widehat{\mu}_{N}^{k}(t)\}$ and $i\neq i^*.$ If the agent $k$ chooses option $i^*$ at time step $t+1$ we have $Q_{i^*}^k(t)> Q_{i}^k(t).$ Thus we  have $\widehat{\mu}_{i}^{k}(t)>\widehat{\mu}_{i^*}^{k}(t)$ and $C_{i}^k(t)<C_{i^*}^k(t).$ 

 Note that for some $\beta_i^k(t)>0$ we have
\begin{align*}
&\mathbb{P}\left(\widehat{\mu}_{i^*}^{k}(t)\geq \mu_{i^*}-C_{i^*}^{k}(t),\widehat{\mu}_{i}^k(t)\geq \widehat{\mu}_{i^*}^k(t)\right)=\beta_i^k(t)\\&+\mathbb{P}\left(\widehat{\mu}_{i^*}^{k}(t)\geq \mu_{i^*}-C_{i^*}^{k}(t),\widehat{\mu}_{i}^k(t)\geq \widehat{\mu}_{i^*}^k(t),N_{i^*}^k(t)\geq \beta_i^k(t)\right).
\end{align*}

Let $\beta_{i}^k(t)=\frac{8\sigma_{i}(\xi+1)}{\Delta^2_{i}}\log t.$
Then we have
\begin{align*}
&\sum_{t=1}^T\mathbb{P}\left(\widehat{\mu}_{i^*}^{k}(t)\geq \mu_{i^*}-C_{i^*}^{k}(t),\widehat{\mu}_{i}^k(t)\geq \widehat{\mu}_{i^*}^k(t)\right)=\beta_i^k(T)\\&+\sum_{i=1}^T\mathbb{P}\left(\widehat{\mu}_{i^*}^{k}(t)\geq \mu_{i^*}-C_{i^*}^{k}(t),\widehat{\mu}_{i}^k(t)\geq \widehat{\mu}_{i^*}^k(t),N_{i^*}^k(t)\geq \beta_i^k(t)\right).
\end{align*}
Since $C_i^k(t)<C_{i^*}^k(t)$ we have
\begin{align*}
&\sum_{i=1}^T\mathbb{P}\left(\widehat{\mu}_{i^*}^{k}(t)\geq \mu_{i^*}-C_{i^*}^{k}(t),\widehat{\mu}_{i}^k(t)\geq \widehat{\mu}_{i^*}^k(t),N_{i^*}^k(t)\geq \beta_i^k(t)\right)\\
&\leq 
\sum_{t=1}^T\mathbb{P}\left(\widehat{\mu}_i^k(t)\geq \mu_i+C_i^k(t)\right)\\
&\leq \beta_i^k(T)+\frac{1}{2^{\xi}\log \zeta}\left(\frac{\log (d_k+1)}{\xi}+\frac{2}{\xi-1}\right)\\
&+\frac{1}{\log \zeta}(1+\log (d_k+1))\\
&+\frac{1}{T^{\xi-1}\log \zeta}\left(\frac{\log (d_k+1)}{T\xi}+\frac{1}{\xi-1}\right).
\end{align*}

Then we have
\begin{align*}
&\sum_{i=1}^T\mathbb{P}\left(\widehat{\mu}_{i^*}^{k}(t)\geq \mu_{i^*}-C_{i^*}^{k}(t),\exists i, s.t. (\widehat{\mu}_{i}^k(t)\geq \widehat{\mu}_{i^*}^{k}(t)\right)\leq\\
&\sum_{i=1}^N\frac{8\sigma_{i}(\xi+1)}{\Delta^2_{i}}\log T +\frac{N-1}{2^{\xi}\log \zeta}\left(\frac{\log (d_k+1)}{\xi}+\frac{2}{\xi-1}\right)\\
    &+\frac{N-1}{\log \zeta}(1+\log (d_k+1))\\
&+\frac{N-1}{T^{\xi-1}\log \zeta}\left(\frac{\log (d_k+1)}{T\xi}+\frac{1}{\xi-1}\right).
\end{align*}

\end{proofoftheorem}

 
\end{appendix}
 

\bibliographystyle{IEEEtran}
\bibliography{MAMAB}

\newpage
\end{document}